# Global Independence of Irrelevant Alternatives, State-Salient Decision Rules and the Strict Condorcet Choice Function

by

Somdeb Lahiri

ORCID: https://orcid.org/0000-0002-5247-3497

PD Energy University, Gandhinagar (EU-G), India (formerly).

(somdeb.lahiri@gmail.com)

May 12, 2023.

## Abstract

We present a simple proof of a well-known axiomatic characterization of state-salient decision rules, using Weak Dominance Criterion and Global Independence of Irrelevant Alternatives. Subsequently we provide a simple axiomatic characterization of the Strict-Condorcet choice function on the domain of all preference profiles that have a strict-Condorcet winner, assuming that if the first two ranks are "occupied by the same two alternatives in all states of nature", then the chosen alternative will be the one from these two that is preferred to the other with probability greater than half- provided such an alternative exists. We also show that this result is not valid if we extend the domain to the set of all preference profiles that have a unique weak-Condorcet winner.

## 1. Introduction:

The motivation for the framework of this paper i.e., decision making with state (or criteria) dependent strict rankings of alternatives, can be found in Lahiri (2019). Our mathematical framework is the same as that of mathematical voting theory originating in the work of Pattanaik (1970), although our interpretation (as in Lahiri (2019)) is different.

To begin with, we provide a simple proof of an axiomatic characterizations of state-salient decision rules, which in the context of the mathematically identical framework of voting theory, due to Pattanaik (1970) and then Denicolo (1985), is about a voting rule being dictatorial if it satisfies the Weak Dominance Criterion and Global Independence of Irrelevant Alternatives. Our proof is different from the one available in Denicolo (1985). The first step of the proof is motivated by the "Extremal Lemma" of Geanakopolos (2005). To the best of our knowledge, there is no other direct and simpler proof of the result we are concerned with here, although our proof has similarities with the proof of Theorem A in Reny (2001), which appeals to Monotonicity instead of Global Independence of Irrelevant Alternatives.

Instead of providing the well-known axiomatic characterization of the Condorcet choice function for the case of two alternatives due to May (1952), the excellent short book by W.D. Wallis (i.e., Wallis (2014)), premises the entire presentation on choice theory with state (or criterion) dependent preferences, on the assumption that in the case of choice between two alternatives, if there exists an alternative that is preferred to the other with probability at least half, then the former is chosen. In fact, a major justification for the Global Independence of

Irrelevant Alternatives assumption in Arrowian choice theory, is that it is satisfied by the above-mentioned choice procedure in pairwise comparisons of alternatives, regardless of the number of alternatives on offer. An alternative that is preferred with probability greater than half to all other alternatives in pair-wise comparisons is said to be a strict Condorcet-winner. By its definition, if such an alternative exists, it must be unique. The problem with making such an assumption, is the one posed by the well Condorcet paradox, which shows that there are preference profiles with three or more alternatives, where such an alternative does not exist. We therefore restrict our analysis to those preference profiles for which such alternatives exist and provide a brief and easily understandable axiomatic characterization of the resolute (singleton valued) choice function which selects for each such profile its unique strict Condorcet-winner. The resulting choice function is called the strict Condorcet choice function.

An alternative is said to be Pareto-dominated if there exists another alternative that is strictly preferred to the former in all states of nature. An alternative that is not Pareto-dominated, is said to be a Pareto-undominated alternative. Clearly a strict-Condorcet winner is a Pareto-undominated alternative.

There is a considerably more complicated axiomatic characterization of the strict-Condorcet choice function in Campbell and Kelly (2003), whose significance and worth are not being questioned here. Like Campbell and Kelly (2003), in our axiomatic characterization, we refrain from assuming that the resolute choice function selects undominated alternatives. The purpose of our axiomatic characterization is the maximum "simplicity" that is possible, if we assume that if the first two ranks are "occupied by the same two alternatives in all states of nature", then the chosen alternative will be the one from these two that is preferred to the other with probability greater than half- provided such an alternative exists. We call this property, "Most Probable at the Top". The other property we invoke in our axiomatic characterization is the Global Independence of Irrelevant Alternatives used in Denicolo (1985). A detailed investigation of the several Independence of Irrelevant Alternatives axioms, preceding the one due to Denicolo (1985) can be found in Ray (1973). A more recent pre-occupation with one of the versions of the axiom discussed in Ray (1973), is the work of Salles (2023).

A weak-Condorcet winner is an alternative that is preferred to all other alternatives with probability at least half and hence a preference profile that has a weak- Condorcet winner may have more than one such alternative. We provide an example with four alternatives and four states of nature, to show that the choice function that selects weak-Condorcet winners on the set of all preference profiles that have a non-empty set of weak-Condorcet winners, does not satisfy Global Independence of Irrelevant Alternatives. On this larger domain but with six states of nature, there is no resolute choice function that always selects a weak-Condorcet winner and satisfies Global Independence of Irrelevant Alternatives. In fact, on the domain of preference profiles with six states of nature and at least three alternatives where each preference profile has a unique weak Condorcet winner, the resolute choice function that always selects the unique weak-Condorcet winner, violates the Global Independence of Irrelevant Alternatives, thereby implying the previous result as well as Global Independence of Irrelevant Alternatives are very restrictive.

After obtaining the results mentioned above, we introduce two conditions on a general domain. The first says that if a strict-Condorcet winner exists, it should be the only alternative to be chosen. The second says that if a strict-Condorcet loser (i.e., an alternative that is ranked worse with probability greater than half than all other alternatives in pair-wise comparisons) exists, then it should not be chosen. Clearly the strict-Condorcet choice function satisfies both conditions. The Plurality rule satisfies neither and the rule that selects the Borda-winner, does not satisfy the former but satisfies the latter. We provide two examples to show the claimed violations for Plurality and Borda. Merlin (2003) presents a survey of the literature on these and other voting rules.

In a final section of this paper, we extend the definition of Strict-Condorcet choice function to the situation where the number of states of nature is variable. An axiomatic characterization of this choice function is an immediate consequence of our axiomatic characterization of the Strict-Condorcet choice function for a fixed set of states of nature.

## 2. The Framework:

Let X be a non-empty finite set of alternatives containing at least three alternatives, from which a decision maker is required to choose exactly one alternative. Unless otherwise required, we will represent the cardinality of X by 'm'. There is a non-empty finite set of future states of nature, which unravel only after the decision maker has made his/her choice. Let N(n) be the set of first 'n' positive integers for some positive integer "greater than or equal to 2. N(n) denotes the set of **states of nature**, whose cardinality is equal to n.

A **strict (preference) ranking** on X is a reflexive, complete/connected/total, transitive and anti-symmetric binary relation on X. Generally, a preference relation is denoted by R. If for $x,y \in X$, it is the case that $(x,y) \in R$, then we shall denote it by xRy. If $x,y \in X$ with $x \neq y$, then xRy is interpreted as x is **strictly preferred** to y.

Let $\mathcal{L}$ denote the set of all strict rankings on X.

Given $R \in \mathcal{L}$ and $x \in X$, let rk(x,R) = cardinality of $\{y \in X | yRx\}$. rk(x,R) is said to be the **rank** of x at R.

Thus, a (strict) ranking R of the alternatives in X can be viewed as a one-to-one function from X to N(n) whose image at each $x \in X$ is said to be the rank of x at the strict ranking and for $x,y \in X$ with $x \neq y$, the rank of x is less than the rank of y if and only if x is preferred to y at the strict ranking.

A **(state-dependent) preference profile** denoted $R_{N(n)}$ is a function from N(n) to $\mathcal{L}$. $R_{N(n)}$ is represented as the array $<R_i | i \in N(n)>$, where $R_i$ is the strict ranking of the decision maker of the alternatives in X in state of nature i. The **set of all preference profiles** is denoted $\mathcal{L}^{N(n)}$.

A (strict) preference profile $R_{N(n)}$ can be represented as an m×n matrix whose $j^{th}$ column for $j \in \{1,\ldots,n\}$ is the strict ranking realized in state of nature 'j', i.e., that entry at the intersection of the $i^{th}$ row and $j^{th}$ column of $R_{N(n)}$, called the $(i,j)^{th}$ entry of $R_{N(n)}$, is the alternative that is ranked $i^{th}$ in state of nature 'j'.

A domain $\mathcal{D}^{N(n)}$ is any non-empty subset of $\mathcal{L}^{N(n)}$.

Let $\Psi(X)$ denote the set of all non-empty subsets of X.

A **Choice Function** (CF) on (a domain) $\mathcal{D}^{N(n)}$ is a function f from $\mathcal{D}^{N(n)}$ to $\Psi(X)$.

In related literature, a choice function is often referred to as a "**winner selection rule**" and the set assigned by the rule to a preference profile as the **choice set** or the **winner set** at the profile.

A CF f on $\mathcal{D}^{N(n)}$ is said to satisfy the **Weak Dominance Criterion** (WDC) if for all x,y∈X with x≠y and $R_{N(n)} \in \mathcal{D}^{N(n)}$: [$xR_iy$ for all i∈N(n)] implies [$y \notin f(R_{N(n)})$].

The above criterion is available and discussed in Lahiri (2019).

The content of the following assumption is available in Denicolo (1985).

A CF f on $\mathcal{D}^{N(n)}$ is said to satisfy **Global Independence of Irrelevant Alternatives** (GIIA) if for all $R_{N(n)}, R'_{N(n)} \in \mathcal{D}^{N(n)}$ and x,y∈X with x ≠ y: [$R_i|\{x,y\} = R'_i|\{x,y\}$ for all i∈N(n), $x \in f(R_{N(n)}), y \notin f(R_{N(n)})$] implies [$y \notin f(R'_{N(n)})$].

The following definition of a choice function on the domain $\mathcal{L}^{N(n)}$ can be found in Lahiri (2019).

A CF f on $\mathcal{L}^{N(n)}$ is said to be a **state-salient decision rule** (S-SDR), if there exists a state of nature i, such that for all x,y∈X with x ≠ y and $R_{N(n)} \in \mathcal{L}^{N(n)}$: $xR_iy$ implies $y \notin f(R_{N(n)})$.

It follows from the definition of an S-SDR that the state of nature that decides the outcome is unique.

Given a preference profile $R_{N(n)}$, an alternative x is a **strict Condorcet-winner** at $R_{N(n)}$ if for all y∈X\{x}, the cardinality of {i∈N(n)| $xR_iy$}> the cardinality of {i∈N(n)| $yR_ix$}.

It is easy to see that if a strict Condorcet-winner exists at $R_{N(n)}$, then it must be unique.

Given a preference profile $R_{N(n)}$, it is quite possible that a strict Condorcet-winner does not exist at $R_{N(n)}$.

Let $\mathcal{D}^{N(n)}_{Strict\ Condorcet} = \{R_{N(n)} \in \mathcal{L}^{N(n)}|$ there is a strict Condorcet-winner at $R_{N(n)}\}$.

For instance, any profile in which some x∈X is ranked first in more than $\frac{n}{2}$ states of nature, belongs to $\mathcal{D}^{N(n)}_{Strict\ Condorcet}$. It is not necessary that 'n' is an odd integer.

A **Resolute Choice Function** (RCF) on a domain $\mathcal{D}^{N(n)}$ is a "<u>singleton</u>-valued" choice function on $\mathcal{D}^{N(n)}$.

If f is an RCF on $\mathcal{D}^{N(n)}$, then for $R_{N(n)} \in \mathcal{D}^{N(n)}$ instead of writing $f(R_{N(n)}) = \{x\}$, we will write $f(R_{N(n)}) = x$.

The **Strict-Condorcet Rule** denoted S-C on $\mathcal{D}^{N(n)}_{Strict\ Condorcet}$ is the RCF that for all $R_{N(n)} \in \mathcal{D}^{N(n)}_{Strict\ Condorcet}$, S-C($R_{N(n)}$) is the strict Condorcet-winner at $R_{N(n)}$.

**3. An axiomatic characterization of State Salient Decision Rules:**

For x∈X, let $\mathfrak{D}^{N(n)}(x) = \{R_{N(n)} \in \mathcal{L}^{N(n)}|$ for all i∈N(n), $rk(x,R_i) \in \{1,m\}\}$.

**Theorem 1:** A CF on $\mathcal{L}^{N(n)}$ satisfies WDC and GIIA if and only if it is an S-SDR.

**Proof:** Let f be a CF on $\mathcal{L}^{N(n)}$.

If f is an S-SDR, then it is easy to see that it satisfies WDC and GIIA.

Hence suppose f satisfies WDC and GIIA.

**Step 1:** Let $x \in X$.

Since by WDC, $rk(x,R_i) = 1$ for all $i \in N(n)$ implies $f(R_{N(n)}) = \{x\}$, let $j(x) = \min\{j \in N(n)|$ there exists $R_{N(n)} \in \mathfrak{D}^{N(n)}(x)$ satisfying $\{i|rk(x,R_i) = 1\} = \{i \in N(n)|i \leq j\}$ and $f(R_{N(n)}) = \{x\}\}$.

By GIIA, if for $R_{N(n)} \in \mathcal{L}^{N(n)}$ and $w \in X\setminus\{x\}$, we have $\{i \in N(n)|xR_iw\} = \{i \in N(n)| i \leq j(x)\}$, then, $w \notin f(R_{N(n)})$.

Let $y,z \in X\setminus\{x\}$ with $y \neq z$ and $R_{N(n)} \in \mathcal{L}^{N(n)}$ such that $\{i|xR_iz\} = \{i \in N(n)| i \leq j(x)\}$, $zR_iy$ for all $i \in N(n)$ and $xR_jy$ for all $i \leq j(x)$, the relationship between $x$ and $y$ being otherwise arbitrary. Thus, $xR_izR_iy$ for all $i \leq j(x)$ and $zR_ix$ for all $i > j(x)$ if $j(x) < n$. Further, suppose $rk(x,R_i)$, $rk(y,R_i)$, $rk(z,R_i) \in \{1,2,3\}$ for all $i \in N$.

By WDC, $f(R_{N(n)}) \subset \{x,y,z\}$ and by WDC (once again) $y \notin f(R_{N(n)})$.

By GIIA, $z \notin f(R_{N(n)})$.

Thus $f(R_{N(n)}) = \{x\}$.

By GIIA, for all $R_{N(n)} \in \mathcal{L}^{N(n)}$ and $y \in X\setminus\{x\}$, $xR_jy$ for all $i \leq j(x)$ implies $y \notin f(R_{N(n)})$.

Thus, for all $x \in X$, there exists $j(x) \in N(n)$ such that for all $R_{N(n)} \in \mathcal{L}^{N(n)}$ and $y \in X\setminus\{x\}$, $xR_jy$ for all $i \leq j(x)$ implies $y \notin f(R_{N(n)})$.

**Step 2:** Let $x \in X$ and $y,z \in X\setminus\{x\}$ with $y \neq z$ and let $R_{N(n)} \in \mathcal{L}^{N(n)}$ be such that $xR_iz$, $xR_iy$, relationship between $z$ and $y$ arbitrary for all $i < j(x)$ (if $j(x) > 1$), $zR_{j(x)}xR_{j(x)}y$, $zR_ix$, for all $i > j(x)$, relationship between $x, y$ and $y, z$ arbitrary for $i > j(x)$. Further, suppose $rk(x,R_i)$, $rk(y,R_i)$, $rk(z,R_i) \in \{1,2,3\}$ for all $i \in N(n)$.

By WDC, $f(R_{N(n)}) \subset \{x,y,z\}$.

By the definition of $j(x)$, and GIIA combined with $xR_iz$ for all $i < j(x)$, $zR_ix$ for all $i \geq j(x)$, we get that $f(R_{N(n)}) \neq \{x\}$.

Since, $xR_iy$ for all $i \leq j(x)$ by the concluding statement of Step 1, we get that $y \notin f(R_{N(n)})$.

Thus, $z \in f(R_{N(n)})$.

By GIIA, for all $y,z \in X\setminus\{x\}$ with $y \neq z$: $[R_{N(n)} \in \mathcal{L}^{N(n)}$ and $zR_{j(x)}y]$ implies $[y \notin f(R_{N(n)})]$.

**Step 3:** We have shown the for each $x \in X$, there exists $j(x) \in N$ such for all $y,z \in X\setminus\{x\}$ with $y \neq z$: $[R''_{N(n)} \in \mathcal{L}^{N(n)}$ and $zR''_{j(x)}y]$ implies $[y \notin f(R''_{N(n)})]$.

Towards a contradiction suppose that for some $x,y \in X$ with $x \neq y$, $j(x) \neq j(y)$.

Let $R''_{N(n)}$ be any profile satisfying $yR''_{j(x)}w$ for all $w \in X\setminus\{y\}$ and $xR''_{j(y)}w$ for all $w \in X\setminus\{x\}$.

Since, $yR''_{j(x)}w$ for all $w \in X\setminus\{y\}$, we get $f(R''_{N(n)}) \cap (X\setminus\{y\}) = \phi$, whence $f(R''_{N(n)}) = \{y\}$.

At the same time, since $xR''_{j(y)}w$ for all $w \in X\setminus\{x\}$, we get $f(R''_{N(n)}) \cap (X\setminus\{x\}) = \phi$, whence $f(R''_{N(n)}) = \{x\}$.

Thus, $x = y$, leading to a contradiction.

Thus, there exists $j \in N(n)$ such that $j(x) = j$ for all $x \in X$, i.e., f is an S-SDR. Q.E.D.

The following example concerns a CF different from S-SDR that satisfies GIIA. However, it does not satisfy WDC.

**Example 1:** Given $R_{N(n)} \in \mathcal{L}^{N(n)}$, let $f(R_N)$ consist of the unique alternative that is ranked last in state of nature 1. It is easy to see that f violates WDC and Monotonicity as defined in Reny (2001), but satisfies GIIA.

An example of a CF-in fact an RCF- that does not satisfy GIIA is the following.

**Example 2:** Let m = 3 with $X = \{x_1, x_2, x_3\}$ and n = 5.

$$\text{Let } R^1_{N(5)} = \begin{bmatrix} x_1 & x_1 & x_2 & x_3 & x_3 \\ x_2 & x_2 & x_3 & x_2 & x_2 \\ x_3 & x_3 & x_1 & x_1 & x_1 \end{bmatrix} \text{ and } R^2_{N(5)} = \begin{bmatrix} x_1 & x_1 & x_2 & x_2 & x_2 \\ x_2 & x_2 & x_3 & x_3 & x_3 \\ x_3 & x_3 & x_1 & x_1 & x_1 \end{bmatrix}.$$

Let $\mathcal{D}^{N(5)}$ be any domain containing $R^1_{N(5)}$ and $R^2_{N(5)}$ and let f be the RCF on $\mathcal{D}^{N(5)}$ such that for all $R_{N(5)} \in \mathcal{D}^{N(5)}$, $f(R_{N(5)}) = x_i$ if and only if $\{i\} =$
$\underset{h \in \{1,2,3\}}{\text{argmin}}\{x_h \in \underset{x_k \in X}{\text{argmax}} \; cardinality\{j \in N(5) | rk(x_k, R_j) = 1\}$, i.e. $f(R_{N(5)})$ is the Plurality winner with the least sub-script.

Thus, $f(R^1_{N(5)}) = x_1$ and $f(R^2_{N(5)}) = x_2$, although $R^1_j|\{x_1,x_2\} = R^2_j|\{x_1, x_2\}$ for all $j \in \{1,\ldots,5\}$. Thus, f violates GIIA.

Note that both $R^1_{N(5)}$ and $R^2_{N(5)}$ belong to $\mathcal{D}^{N(5)}_{Strict\;Condorcet}$ with $x_2$ being the strict- Condorcet winner at both profiles.

A CF f on $\mathcal{D}^{N(n)}$ is said to satisfy **Weak Monotonicity** if for all $R_N, R'_N \in \mathcal{D}^{N(n)}$ and $x \in f(R_N)$: $[\{i|xR_iy\} \subset \{i|xR'_iy\}$ for all $y \in X \setminus \{x\}]$ implies $[x \in f(R'_N)]$.

The statement $[\{i|xR_iy\} \subset \{i|xR'_iy\}$ for all $y \in X \setminus \{x\}]$ is equivalent to the statement $[\{y \in X \setminus \{x\}|xR_iy\} \subset \{y \in X \setminus \{x\}|xR'_iy\}$ for all $i \in N(n)]$.

A CF f on $\mathcal{D}^{N(n)}$ is said to be **Resolute for Pairs** if for all $R_N \in \mathcal{D}^{N(n)}$ and $x,y \in X$ with $x \neq y$: $[rk(x,R_i), rk(y,R_i) \in \{1,2\}$ for all $i \in N]$ implies $[f(R_N) = \{x\}$ or $f(R_N) = \{y\}]$.

An S-SDR is Resolute for Pairs and satisfies Weak Monotonicity.

**Observation:** Suppose f on $\mathcal{L}^{N(n)}$ is Resolute for Pairs and satisfies Weak Monotonicity. Then f satisfies GIIA.

Towards a contradiction suppose that for some $R_{N(n)}, R'_{N(n)} \in \mathcal{L}^{N(n)}$ and $x,y \in X$ with $x \neq y$ it is the case that $R_i|\{x,y\} = R'_i|\{x,y\}$ for all $i \in N(n)$, $x \in f(R_{N(n)})$, $y \notin f(R_{N(n)})$ and yet $y \in f(R'_{N(n)})]$.

Let $R''_{N(n)} \in \mathcal{L}^{N(n)}$ be such that $R''_i|\{x,y\} = R_i|\{x,y\}$, $R''_i|(X \setminus \{x,y\}) = R_i|(X \setminus \{x,y\})$, $rk(x, R''_i)$, $rk(y, R''_i) \in \{1,2\}$ for all $i \in N(n)$.

Thus, $rk(x, R''_i), rk(y, R''_i) \in \{1,2\}$ for all $i \in N(n)$.

Since f is Resolute for Pairs, either $f(R''_{N(n)}) = \{x\}$ or $f(R''_{N(n)}) = \{y\}$

Since $\{i|xR_iw\} \subset \{i|xR_i''w\}$ for all $w \in X\setminus\{x\}$, by Weak Monotonicity and Resolute for Pairs, we get $f(R''_{N(n)}) = \{x\}$.

Since $\{i|yR'_iw\} \subset \{i|yR_i''w\}$ for all $w \in X\setminus\{y\}$, by Weak Monotonicity and Resolute for Pairs we get $f(R''_{N(n)}) = \{y\} \neq \{x\}$, leading to a contradiction.

Thus, f satisfies GIIA.

Thus, and in view of theorem 1, we get the following result.

**Corollary 1 of Theorem 1:** A CF on $\mathcal{L}^{N(n)}$ satisfies WDC, Weak Monotonicity and is Resolute for Pairs if and only if it is an S-SDR.

**Note:** The corollary and the discussion preceding it, is not valid in the absence of the Resolute for Pairs assumption. For, let f be the CF on $\mathcal{L}^{N(n)}$ such that for all $R_{N(n)} \in \mathcal{L}^{N(n)}$, $f(R_{N(n)}) = \{x \in X: rk(x,R_i) = 1 \text{ for some } i \in N(n)\}$. It is easy to see that f satisfies Weak Monotonicity but not Resolute for Pairs. For the latter consider the set of alternatives $X = \{x,y,z\}$, $n = 3$ and the preference profile $\begin{bmatrix} x & x & y \\ y & y & x \\ z & z & z \end{bmatrix}$, the chosen set of alternatives for which is $\{x,y\}$, even though the first two ranks are occupied by either x or y in every state of nature. It does not satisfy GIIA either, as is apparent in the situation with $X = \{x,y,z\}$ and $n = 3$ represented below.

Let $R_{N(3)} = \begin{bmatrix} x & x & y \\ y & y & z \\ z & z & x \end{bmatrix}$. $f(R_{N(3)}) = \{x,y\}$ and $z \notin f(R_{N(3)})$.

Now let $R'_{N(3)} = \begin{bmatrix} x & x & z \\ y & y & y \\ z & z & x \end{bmatrix}$. Clearly $R'_i|\{x,z\} = R_i|\{x,z\}$ for all $i \in N(3)$ and yet $z \in f(R'_{N(3)}) = \{x,z\}$ leading to a violation of GIIA.

The following property is defined in chapter 3 of Taylor (2005).

A CF f on $\mathcal{D}^{N(n)}$ is said to satisfy **Down Monotonicity** if for all $R_N, R'_N \in \mathcal{D}^{N(n)}$ $i \in N(n)$, $x \in f(R_N)$ and $y \neq x$: $[R_j = R'_j$ for all $j \neq i$, $R_i|(X\setminus\{y\}) = R'_i|(X\setminus\{y\})$, $rk(y, R'_i) = rk(y,R_i) + 1]$ implies $[x \in f(R'_N)]$.

It is easy to see that an S-SDR satisfies Down Monotonicity.

The proof of the following proposition is provided in an Appendix to this paper.

**Proposition 1:** On $\mathcal{L}^{N(n)}$ Weak Monotonicity is equivalent to Down Monotonicity.

In view of this proposition and Corollary 1 of Theorem 1 we get the following result.

**Corollary 2 of Theorem 1:** A CF on $\mathcal{L}^{N(n)}$ satisfies WDC, Down Monotonicity and is Resolute for Pairs if and only if it is an S-SDR.

**4. Axiomatic Characterization of the Strict Condorcet Choice Rule:**

From this section onwards we assume that "all states of nature are **equiprobable**".

The following is a property, required for the main result of this section, which we refer to as "Most Probable at the Top".

A CF f on $\mathcal{D}^{N(n)}$ is said to satisfy **Most Probable at the Top (MPT)** if for all $R_{N(n)} \in \mathcal{D}^{N(n)}$ and x,y∈X with x ≠ y such that $xR_iw$, $yR_iw$ for all i∈N(n) and w∈X\{x,y}: [cardinality of {i∈N(n)| $xR_iy$}> cardinality of {i∈N(n)| $yR_ix$}] implies [f($R_N$) = {x}].

Recall that if f is an RCF, then instead of writing f as a singleton set, we have agreed to represent f as a function whose co-domain is X, i.e., instead of writing $f(R_{N(n)})$ = {x}, we have agreed to write $f(R_{N(n)})$ = x

**Theorem 2**: A CF f with domain $\mathcal{D}^{N(n)}_{Strict\ Condorcet}$ satisfies GIIA and MPT if and only if f = S-C.

**Proof:** It is easy to see that on $\mathcal{D}^{N(n)}_{Strict\ Condorcet}$ S-C satisfies MPT.

Suppose, $R_{N(n)} \in \mathcal{D}^{N(n)}_{Strict\ Condorcet}$, with S-C($R_{N(n)}$) = x∈X and let y∈X\{x}. Thus, x is preferred to y with probability greater than $\frac{1}{2}$ at $R_{N(n)}$. Let $R'_{N(n)} \in \mathcal{D}^{N(n)}_{Strict\ Condorcet}$ and suppose $R'_i|\{x,y\} = R_i|\{x,y\}$ for all i∈N(n). Thus, x is preferred to y with probability greater than $\frac{1}{2}$ at $R'_{N(n)}$. Thus, S-C($R'_{N(n)}$)≠ y. Thus, S-C satisfies GIIA.

Now suppose, f is an CF on $\mathcal{D}^{N(n)}_{Strict\ Condorcet}$ that satisfies GIIA and MPT.

Let $R_{N(n)} \in \mathcal{D}^{N(n)}_{Strict\ Condorcet}$ with strict Condorcet-winner being x and y∈f($R_{N(n)}$). Towards a contradiction, suppose that y ≠ x.

Let $R^{\{x,y\}}_{N(n)}$ be the preference profile such that:

(i) For all i∈N(n), $R^{\{x,y\}}_i|\{x,y\} = R_i|\{x,y\}$ and $R^{\{x,y\}}_i|(X\setminus\{x,y\}) = R_i|(X\setminus\{x,y\})$;

(ii) For all i∈N(n): rk(x, $R^{\{x,y\}}_i$)∈{1,2} and rk(y, $R^{\{x,y\}}_i$)∈{1,2}.

It is easy to see that $R^{\{x,y\}}_{N(n)} \in \mathcal{D}^{N(n)}_{Strict\ Condorcet}$, with strict Condorcet winner being x.

Since, $xR^{\{x,y\}}_iw$, $yR^{\{x,y\}}_iw$ for all i∈N(n) and w∈X\{x,y}: [cardinality of {i∈N(n)| $xR^{\{x,y\}}_iy$}> cardinality of {i∈N(n)| $yR^{\{x,y\}}_ix$}] along with MPT implies $f(R^{\{x,y\}}_{N(n)})$ = {x}.

However, x∈$f(R^{\{x,y\}}_{N(n)})$, y∉$f(R^{\{x,y\}}_{N(n)})$, $R^{\{x,y\}}_i|\{x,y\} = R_i|\{x,y\}$ for all i∈N(n) implies by GIIA that y∉f($R_{N(n)}$) leading to a contradiction.

This proves the theorem. Q.E.D.

Given $R_{N(n)} \in \mathcal{D}^{N(n)}$, x∈X is said to be a **weak-Condorcet winner** at $R_{N(n)}$ if for all y∈X\{x}: cardinality of {j∈N(n)| rk(x,$R_j$) < rk(y,$R_j$)} ≥ $\frac{n}{2}$.

Let $\mathcal{D}^{N(n)}_{Weak\ Condorcet}$ = {$R_{N(n)} \in \mathcal{L}^{N(n)}$|there exists a weak Condorcet winner at $R_{N(n)}$}.

Let C be the CF on $\mathcal{D}_{Weak\ Condorcet}^{N(n)}$ such that for all $R_{N(n)} \in \mathcal{D}_{Weak\ Condorcet}^{N(n)}$, $C(R_{N(n)})$ is the set of weak-Condorcet winners at $R_{N(n)}$.

The following example shows that C does not satisfy GIIA.

**Example 3:** Let m = 4 and n = 4.

Let $R_{N(4)} = \begin{bmatrix} x & x & y & z \\ y & y & w & w \\ z & z & x & x \\ w & w & z & y \end{bmatrix}$.

x is preferred to y and z individually with probability $\frac{3}{4}$ and x is preferred to w with probability $\frac{1}{2}$. y and z are individually preferred to w with probability $\frac{3}{4}$. Thus, $C(R_{N(4)}) = \{x\}$.

Let $R'_{N(4)} = \begin{bmatrix} x & x & w & w \\ y & y & y & z \\ z & z & x & x \\ w & w & z & y \end{bmatrix}$.

$R'_{N(4)}$ is obtained from $R_{N(4)}$ by interchanging the positions of the first two alternatives in the third and fourth states of nature.

It is easy to see that $C(R'_{N(4)}) = \{x,w\}$, since w is preferred to x,y,z individually with probability equal to $\frac{1}{2}$, the other probabilities being the same as before.

$C(R_{N(4)}) = \{x\}$, $R_i|\{x,w\} = R'_i|\{x,w\}$ for all $i \in N$ and yet $w \in C(R'_{N(4)}) = \{x,w\}$. Hence C violates GIIA.

Theorem 2 is quite "tight", in that it cannot be extended to $\mathcal{D}_{Weak\ Condorcet}^{N(4)}$, i.e., there does not exist any RCF that always selects a weak-Condorcet winner and satisfies GIIA. In fact, this observation is implied by a much stronger result for which we require the following domain.

Let $\mathcal{D}_{Weak\ Condorcet}^{N(n)*} = \{R_{N(n)} \in \mathcal{D}_{Weak\ Condorcet}^{N(n)} |$ there exists a <u>unique</u> weak Condorcet winner at $R_{N(n)}\}$.

**Proposition 2:** For m = 3, and n = 6, let f be the RCF on $\mathcal{D}_{Weak\ Condorcet}^{N(6)*}$ such that for each $R_{N(6)} \in \mathcal{D}_{Weak\ Condorcet}^{N(6)*}$, $f(R_{N(6)})$ is the unique weak Condorcet winner at $R_{N(6)}$. Then, f violates GIIA.

**Proof:** Let f be the RCF on $\mathcal{D}_{Weak\ Condorcet}^{N(6)*}$ defined in the statement of Proposition 2.

Let $R_{N(6)}^{(1)} = \begin{bmatrix} x & x & x & z & y & y \\ z & z & z & y & z & x \\ y & y & y & x & x & z \end{bmatrix}$.

The unique weak Condorcet winner at $R_{N(6)}^{(1)}$ is 'x' and so $f(R_{N(6)}^{(1)}) = x$.

Let $R_{N(6)}^{(2)} = \begin{bmatrix} x & x & z & z & y & y \\ y & y & x & y & z & z \\ z & z & y & x & x & x \end{bmatrix}$.

Clearly, $R_{N(6)}^{(2)}|\{x,y\} = R_{N(6)}^{(1)}|\{x,y\}$. However the, weak Condorcet winner at $R_{N(6)}^{(2)}$ is 'y' and so $f(R_{N(6)}^{(2)}) = y$, leading to a violation of GIIA. Q.E.D.

Proposition 2 shows that GIIA is extremely restrictive. This perception is reinforced if we consider the next example for which we require the following domain.

$$\mathcal{D}_{Plurality}^{N(n)*} = \{R_{N(n)} \in \mathcal{L}^{N(n)}| \underset{x \in X}{\mathrm{argmax}}\ cardinality\{i \in N(n)|rk(x, R_i) = 1\}\ \text{is a singleton}\}.$$

(For $R_{N(n)} \in \mathcal{L}^{N(n)}$, $\underset{x \in X}{\mathrm{argmax}}\ cardinality\{i \in N(n)|rk(x, R_i) = 1\}$ is the set of Plurality winners at $R_{N(n)}$. $\mathcal{D}_{Plurality}^{N(n)*}$ is the set of preference profiles with a unique Plurality winner.)

**Example 4:** Let m = 3 and n = 7. For X = {x,y,z}, let f be the CF on $\mathcal{D}_{Plurality}^{N(n)*}$ such the for all $R_{N(7)} \in \mathcal{D}_{Plurality}^{N(n)*}$, $f(R_{N(7)})$ =$\underset{w \in X}{\mathrm{argmax}}\ cardinality\{i \in N(7)|rk(w, R_i) = 1\}$. Let

$$R_{N(7)} = \begin{bmatrix} x & x & x & y & y & z & z \\ y & y & y & z & z & y & y \\ z & z & z & x & x & x & x \end{bmatrix} \text{ and } R'_{N(7)} = \begin{bmatrix} x & x & x & y & y & y & y \\ y & y & y & z & z & z & z \\ z & z & z & x & x & x & x \end{bmatrix}.$$

$f(R_{N(7)}) = \{x\}$, $R_i|\{x,y\} = R'_i|\{x,y\}$ for all $i \in N(7)$ and yet $f(R'_{N(7)}) = \{y\}$, violating GIIA.

## 5. Condorcet winner and Anti-Condorcet loser properties:

Given a strict preference profile $R_{N(n)}$, an alternative $x \in X$ is said to a **strict Condorcet-loser** if for all $y \in X\setminus\{x\}$, the cardinality of $\{i \in N(n)| yR_ix\}$ > the cardinality of $\{i \in N(n)| xR_iy\}$.

Two desirable properties for a choice function are the following.

A CF f on $\mathcal{D}^{N(n)}$ is said to satisfy **the strict Condorcet-winner property** if for all $R_{N(n)} \in \mathcal{D}^{N(n)} \cap \mathcal{D}_{Strict\ Condorcet}^{N(n)}$: $f(R_{N(n)}) = \{$the unique strict Condorcet-winner at $R_{N(n)}\}$.

A CF f on $\mathcal{D}^{N(n)}$ is said to satisfy **the Anti-Condorcet-loser property** if for all $R_{N(n)} \in \mathcal{D}^{N(n)}$: $x \in f(R_{N(n)})$ implies x <u>is not</u> a strict Condorcet-loser.

Clearly S-C on $\mathcal{D}_{Strict\ Condorcet}^{N(n)}$ satisfies both properties.

For $x \in X$ and $R_{N(n)} \in \mathcal{L}^{N(n)}$ the **Borda-score** of x at $R_{N(n)}$, denoted BS(x,$R_{N(n)}$) = $\sum_{j \in N(n)}[m + 1 - rk(x, R_j)] = n(m+1) - \sum_{j \in N(n)} rk(x, R_j)$.

Given $R_{N(n)} \in \mathcal{L}^{N(n)}$ a **Borda-winner at** $R_N$ is an alternative x such that BS(x,$R_{N(n)}$) $\geq$ BS(y,$R_{N(n)}$) for all $y \in X$.

The following example shows that the Plurality rule (i.e., the rule that selects all alternatives that are ranked first with the highest probability) violates both properties.

**Example 5:** Let m = 4 and n = 5. Suppose X = {x,y,z,w}.

Let $R_{N(5)} = \begin{bmatrix} x & x & y & z & w \\ y & y & z & w & z \\ w & z & w & y & y \\ z & w & x & x & x \end{bmatrix}.$

The Plurality rule selects x which is the unique Condorcet-loser. The unique Condorcet winner is 'y'.

In this example the Borda-score of $x = 11$, the Borda-score of $y = 14$, the Borda-score of $z = 13$ and the Borda-score of $w = 12$.

Thus, in this example the unique Condorcet-winner is also the Borda-winner.

**Note:** A Borda-winner cannot be a Condorcet-loser.

Let $R_{N(n)} \in \mathcal{L}^{N(n)}$ and suppose x is a strict Condorcet-loser.

We have $BS(y, R_{N(n)}) = n(m+1) - \sum_{j \in N(n)} rk(y, R_j)$ for all $y \in X$.

Thus, $\sum_{y \in X} BS(y, R_{N(n)}) = nm(m+1) - \sum_{y \in X}[\sum_{j \in N(n)} rk(y, R_j)] = nm(m+1) - \sum_{j \in N(n)}[\sum_{y \in X} rk(y, R_j)]$.

For all $j \in N$, $\sum_{y \in X} rk(y, R_j) = \frac{m(m+1)}{2}$.

Thus, $\sum_{y \in X} BS(y, R_{N(n)}) = nm(m+1) - \frac{nm(m+1)}{2} = \frac{nm(m+1)}{2}$.

Thus, average Borda-score of the alternatives is $\frac{n(m+1)}{2}$.

However, $\sum_{j \in N(n)} rk(y, R_j) = \sum_{j \in N(n)}[\text{cardinality of } \{z \in X \setminus \{y\} | zR_j y\} + 1] = n + \sum_{z \in X \setminus \{y\}} \text{cardninality of } \{j \in N(n) | zR_j y\}$.

If x is a strict Condorcet-loser then, for all $y \in X \setminus \{x\}$: cardinality of $\{j \in N(n) | yR_j x\} > \frac{n}{2}$ so that $BS(x, R_{N(n)}) < n(m+1) - n - \frac{n(m-1)}{2} = n(m+1) - n[1 + \frac{(m-1)}{2}] = n(m+1) - \frac{n(m+1)}{2} = \frac{n(m+1)}{2}$.

Since the average Borda-score is $\frac{n(m+1)}{2}$ and $BS(x, R_{N(n)}) < \frac{n(m+1)}{2}$, there exists $y \in X \setminus \{x\}$ such that $BS(y, R_{N(n)}) > \frac{n(m+1)}{2} > BS(x, R_{N(n)})$ and hence x cannot be a Borda-winner.

However, a Borda-winner may violate the Condorcet-winner property as the following example reveals.

**Example 6:** Let $m = 4$ and $n = 5$.

Let $R_{N(5)} = \begin{bmatrix} x & x & x & z & w \\ y & y & y & y & y \\ w & w & w & w & z \\ z & z & z & x & x \end{bmatrix}$

The unique strict Condorcet-winner is x.

The Borda-score of $x = 14$, the Borda-score of $y = 15$, the Borda-score of $z = 9$ and the Borda-score of $w = 12$.

Thus, the Borda-winner is y.

**6. Choice functions on domains with variable sets of states of nature:**

In this section, we extend the definition of a Strict-Condorcet Choice function to the situation where the number of states of nature is variable.

Let $\mathfrak{L} = \bigcup_{n \geq 2} \mathcal{L}^{N(n)}$ be the set of all state-dependent preference profiles with arbitrarily many-though a finite number and greater than or equal to two- states of nature.

A **domain** (in this context) is any non-empty subset $\mathfrak{D}$ of $\mathfrak{L}$.

An **Anonymous Choice Function (ACF)** on a domain $\mathfrak{D}$ is a function f: $\mathfrak{D} \to \Psi(X)$ such that for all $R_{N(n)}, R'_{N(r)} \in \mathfrak{D}$: [ $\frac{cardinality\ of\ \{j|R_j=R\}}{n} = \frac{cardinality\ of\ \{j|R'_j=R\}}{r}$ for all $R \in \mathcal{L}$] implies [f($R_{N(n)}$) = f($R'_{N(r)}$)].

The typical way of representing $R_{N(n)} \in \mathfrak{L}$ in such a context is by using a matrix having m rows and the number of columns of the matrix being equal to the cardinality of $\{R \in \mathcal{L} | \{j|R_j = R\} \neq \phi\}$.

If $\{R \in \mathcal{L} | \{j|R_j = R\} \neq \phi\} = \{R^{(1)},\ldots,R^{(K)}\} \subset \mathcal{L}$, for some positive integer K, then clearly $\sum_{k=1}^{K}$ cardinality of $\{j|R_j = R^{(k)}\}$ = n and $R_{N(n)}$ is the m×K matrix such that for each $k \in \{1,\ldots,K\}$, the entry in $j^{th}$ row of the $k^{th}$ column is the $j^{th}$ ranked alternative in $R^{(k)}$. The cardinality of $\{j|R_j = R^{(k)}\}$ is usually written on top of the $j^{th}$ column.

Thus if m = 3 and X = {x,y,z} a typical preference profile in $\mathfrak{L}$ could be represented as follows:

$$\begin{matrix} n_1 & \cdots & n_K \\ \begin{bmatrix} x & \cdots & y \\ y & \cdots & z \\ z & \cdots & x \end{bmatrix} \end{matrix}$$

In the above no two columns of the matrix are equal and thus K $\leq \sum_{k=1}^{K} n_k$. Further for each $k \in \{1,\ldots,K\}$, $n_k$ = cardinality$\{j|R_j = R^{(k)}\}$. Thus, for each $k \in \{1,\ldots,K\}$, $\frac{n_k}{\sum_{h=1}^{K} n_h}$ is the probability of the event that the state-dependent preference is $R^{(k)}$.

Let $\mathfrak{D}_{Strict\ Condorcet} = \bigcup_{n \geq 2} \mathfrak{D}_{Strict\ Condorcet}^{N(n)}$

Consider the function $\mathcal{S} - \mathcal{C}: \mathfrak{D}_{Strict\ Condorcet} \to \Psi(X)$, such that for all $R_{N(n)} \in \mathfrak{D}_{Strict\ Condorcet}$, $\mathcal{S} - \mathcal{C}(R_{N(n)}) = \{$strict Condorcet winner at $R_{N(n)}\}$.

It is easy to verify that $\mathcal{S} - \mathcal{C}$ is an ACF- in fact an ARCF (Anonymous "Resolute" Choice Function).

$\mathcal{S} - \mathcal{C}$ is referred to as **the Strict-Condorcet Choice Function with a Variable Number of States of Nature**.

**Acknowledgment:** This paper is based on two earlier papers by the author (a) A Simple Proof of an Axiomatic Characterization of State Salient Decision Rules (available at: https://drive.google.com/file/d/1K5NaBYkn3K7ftXGoEa9DAR5S8S2fJlTo/view), and (b) The Condorcet Choice Function for Equiprobable States of Nature (available at: https://drive.google.com/file/d/1QZUtCh8cYrlvewuIA61_qMq2t5q-mi43/view). A very warm "Thank you" to Franz Deitrich and Joseph Mullat for their very enlightening comments and suggestions.

## Appendix

**Proposition:** On $\mathcal{L}^{N(n)}$ Weak Monotonicity is equivalent to Down Monotonicity.

**Proof:** It is easy to see that Weak Monotonicity implies Down Monotonicity. Hence let us prove the converse.

Suppose f satisfies Down Monotonicity. Let $R_N, R'_N \in \mathcal{L}^{N(n)}$, $x \in f(R_N)$ and [$\{i|xR_iy\} \subset \{i|xR'_iy\}$ for all $y \in X\setminus\{x\}$]. Thus, $rk(x, R'_i) \leq rk(x, R_i)$ for all $i \in N(n)$. We need to show that $x \in f(R'_N)$.

It is enough to show this for the case where there exists $i \in N(n)$ such that $R_j = R'_j$ for all $j \neq i$, since by repeated application of the result in the special case, we obtain the result for the general case.

Case 1: Suppose $rk(x, R'_i) = rk(x, R_i)$.

Thus, $\{y \in X\setminus\{x\} | xR'_i y\} = \{y \in X\setminus\{x\} | xR_i y\}$ and $\{y \in X\setminus\{x\} | yR'_i x\} = \{y \in X\setminus\{x\} | yR_i x\}$

Then, $R'_i$ can be obtained from $R_i$ by repeatedly interchanging adjacent pairs $y,z \in X\setminus\{x\}$ i.e. pairs $y,z \in X\setminus\{x\}$, which are ranked one after another, with the rank of one alternative in the pair decreasing by 1, and the rank of the other increasing by one. Further, either both y and z are preferred to x or x is preferred to both. $x \in f(R'_N)$

By repeated application of Down Monotonicity we get $x \in f(R'_N)$.

Case 2: Suppose $rk(x, R'_i) < rk(x, R_i)$.

Let $R''_N \in \mathcal{L}^{N(n)}$ be such that $R''_j = R_j$ for all $j \neq i$, $R''_i|(X\setminus\{x\}) = R'_i|(X\setminus\{x\})$ but $rk(x, R''_i) = rk(x, R_i)$.

By Case 1, $x \in f(R''_N)$.

Since $R''_i|(X\setminus\{x\}) = R'_i|(X\setminus\{x\})$ but $rk(x, R''_i) = rk(x, R_i) > rk(x, R'_i)$, $R'_i$ is obtained from $R''_i$ by interchanging (if necessary repeatedly) the position of x in the $i^{th}$ state of nature with the alternative immediately above it, the rank of the latter in the $i^{th}$ state of nature thereby increasing by one.

By, if necessary repeated, application of Down Monotonicity, we get $x \in f(R'_N)$.

This proves the proposition. Q.E.D.